# Resource Guide for Teaching Post-Quantum Cryptography

Joshua Holden[a]*

[a]*Department of Mathematics, Rose-Hulman Institute of Technology, Terre Haute, IN, USA*

5500 Wabash Ave., Terre Haute, IN 47803, USA. E-mail: holden@rose-hulman.edu

Abstract:  Public-key cryptography has become a popular way to motivate the teaching of concepts in elementary number theory, abstract algebra, and introduction to proof courses, as well as in cryptography courses. Unfortunately, many experts expect quantum computers to make common forms of public-key cryptography obsolete in the near future.  Fortunately, there are several systems being evaluated to replace RSA and the other systems we currently use.  While some of the systems are too complicated to be good examples in introductory courses, others are either quite manageable or have simplified versions which are manageable.  This article gives a tour of the main types of systems under consideration and the teaching resources available for instructors who want to teach them.

1. **Introduction**

On July 22, 2020, the National Institute of Standards and Technology announced the list of 7 finalists and 8 alternates advancing to the third round of its "post-quantum cryptography standardization process". The modern idea of public-key cryptography is based on "hard problems" --- puzzles which are easy to construct and hard to solve. For example, RSA cryptography is based on the fact that it is easy to multiply two large whole numbers together, and hard to take the result and factor it back into the original two numbers. Hard number theory problems, such as factoring and discrete logarithms, are involved in the vast majority of encrypted messages on the Internet today. Unfortunately, there are signs that quantum computers which can solve these problems quickly may be built in the near future.[1] Cryptographers are looking in all corners of mathematics for more hard problems that quantum computers won't be able to defeat.

There are five major types of hard problems under consideration: lattice problems, code problems, multivariable polynomial problems, hash function problems, and elliptic curve isogenies. Unfortunately for those trying to prioritize, the five types are very different and it's reasonably certain that at least one of them is never going to get standardized. I'll talk first about survey resources and then say a little bit about each type of problem and individual resources.

If you already teach cryptography and want to get yourself up to speed on post-quantum, the standard survey is *Post-Quantum Cryptography* (Bernstein, Buchmann, and Dahmen 2009). It has a chapter on each of most of the main types of post-quantum systems. (Admittedly some are better than others.) It doesn't cover the very latest systems, but most of the ones NIST is evaluating are very similar to ones in the book. I thought the NIST "Status Report on the Second Round of the NIST Post-Quantum Cryptography Standardization Process" (Alagic et al. 2020) was also very helpful for keeping the schemes straight, although it's a little heavy on jargon.

---

[1] Symmetric-key cryptography, as opposed to public-key cryptography, is much less vulnerable to quantum computers. Experts believe that merely using longer keys in the same cryptosystems will be sufficient to keep these systems secure.

As an overview for students, the best options are the 4th edition of *Cryptography: Theory and Practice* (Stinson and Paterson 2018), and the 8th edition of *Cryptography and Network Security: Principles and Practice* (Stallings 2019). Stinson has a brief but complete overview of the rationale for PQC and the five major types under consideration, then four short sections with at least one complete cryptosystem from each of the first four types. Again, it's not always the latest systems but I think they are good choices pedagogically. Stallings has a good overview with excellent motivation but very short on technical details --- it's more suitable for a computer security class than a cryptography class per se.

## 2. Lattice-based Cryptography

Lattice problems are the easiest to describe and NIST has indicated that it's very likely that at least one lattice-based cryptosystem will be standardized. A lattice is a grid of evenly spaced points in a vector space. An example of a problem involving lattices which cryptographers think is hard is the closest vector problem. In this problem you are given a description of the lattice and a point in the space which is not in the lattice. The goal is to find a point in the lattice as close as possible to the given point. In two-dimensional space the solution to this problem has been known since the early $19^{th}$ century, but in, say, five hundred-dimensional space it is still considered very difficult.

The drawback of cryptography that relies only on lattices is that cryptographic keys take up a lot of space, since the description of the space is long.  Therefore, all but one of the systems NIST is considering represent the lattice using algebraic concepts like polynomials and modular arithmetic.  You can still work in 500-dimensional space, but you only need one polynomial to describe the lattice.  That means the secret key to the cipher can be much shorter, which is very convenient.

NTRU is a well-established finalist and old enough to be in standard texts, so it's an obvious place to start with lattice-based systems.  It's covered well in Stinson and even better in *An Introduction to Mathematical Cryptography* (Hoffstein, Pipher, and Silverman 2014) (unsurprisingly, since they invented it), including how to (try to) break it with lattice attacks.  NTRU Prime (a NIST alternate) is newer but is basically a variant on NTRU with some wrapping around it.  FALCON (a finalist) is also related but rather more complicated.

There are several other lattice-based systems in contention:  CRYSTALS-KYBER, CRYSTALS-DILITHIUM, SABER, and FrodoKEM.  These are based on versions of the Learning With Errors (LWE) problem, which is different from the NTRU problem but similar to the closest vector problem.   Stinson covers the original LWE cryptosystem, which is not that different from FrodoKEM.  Another option is *A Course in Cryptography* (Knospe 2019), which has a slightly longer section on LWE, with a few exercises.  It is written from a fairly sophisticated mathematical point of view which might not be suitable for all students.

 The other three lattice-based systems are based on some very new variations of LWE (MLWE and MLWR).  The best overview of LWE (and LWR) that I know of is "A Decade of Lattice Cryptography" (Peikert 2016).  I highly recommend it, but it was last revised in 2016 so it still doesn't have the very latest systems.

3. **Code-based Cryptography**

Another type of cryptosystem which has a good chance of being standardized is based on code problems. Code-based cryptography sounds like a funny name, but that's because the type of codes in question aren't secret codes. They are more specifically known as error-correcting codes, and they were developed to find and correct errors which occur when bits are transmitted across an unreliable communications channel. Cryptography based on codes is actually very similar in some ways to cryptography based on lattices. The distance between two strings of bits is considered to be the number of bits that need to be flipped to get from one to the other. For a random set of codewords it's a hard problem to find the nearest codeword to a given string of bits. So the closest codeword problem, like the closest vector problem, is a hard problem that can be used for cryptography.

One of the code-based NIST finalists, Classic McEliece, is very old by public-key standards. (It is almost as old as RSA, which truly makes it a classic by public-key cryptography standards.) Stinson covers it well. You can find decent expositions in a number of other textbooks, such as *Introduction to Cryptography with Coding Theory* (Trappe and Washington 2020), which generally explains things well but tends to be weak on exercises. (There don't seem to be any specifically on McEliece, although there are some on codes in general.) Knospe goes into more detail, but again is not suitable for all students.

Unfortunately, Classic McEliece has a huge key size. HQC and BIKE (two NIST alternates) are similar but based on quasi-cyclic codes, which are discussed briefly in Bernstein et al. They have smaller key sizes but a much shorter track record.

## 4. Hash-based Cryptography

The next two types of cryptosystems, hash functions and multivariable equations, seem borderline for getting standardized. Both only create digital signatures and they have relatively few candidates still in contention. A hash function is a function which takes an arbitrarily long string as input and produces a short output in a way that is impractically difficult to predict or reverse. The basic idea of using hash functions for digital signatures goes back to the Lamport scheme and Merkle trees which are pretty easy to describe and are both well covered in Stinson and in Bernstein et al. (Pretty well in Stallings, also.) SPHINCS+ (an alternate) is a new system based on those ideas with a lot of tweaks for optimization.[2]

Picnic (an alternate) also uses hash functions, but in the context of non-interactive zero-knowledge proofs of knowledge. A zero-knowledge proof of knowledge is one that allows to prove that you know a certain value, like the output of some function, without revealing what the value is. If you are teaching cryptography, you might already spend time on zero-knowledge proofs and that might make Picnic looking into. I'd say it has a low probability of becoming an important standard, however.

## 5. Multivariable Polynomial Cryptography

Multivariable polynomial cryptography is based on finding values of the variables that make the polynomial equal to zero. With only one variable, mathematicians have pretty much known for centuries how to find solutions or tell if there aren't any. For polynomials with of more than one variable, the problem is much harder. A drawback is again that key sizes for multivariable polynomial ciphers tend to be very large. Also, there doesn't seem to be a unified theory of how to make ciphers out of this hard problem. There are several very different ideas, including "Oil and Vinegar" and "Hidden Field Equations", and it's not clear whether any of them have systemic weaknesses which have not yet been revealed. There are two NIST contenders based on multivariable equations, Rainbow and GeMSS. Rainbow is based on the Oil and Vinegar scheme, and GeMSS is based on Hidden Field Equations. Simple examples of both types are covered well in Stinson. Rainbow itself is covered in Bernstein et al. GeMSS is not, but the closely related HFEv- is.

---

[2] NIST is also moving ahead on standardizing two other variations of the Lamport and Merkle scheme which are "stateful", and therefore harder to use, but are generally agreed to have good speed and security.

6. **Supersingular Elliptic Curve Isogeny-based Cryptography**

The newest type of hard problem, the supersingular isogeny Diffie-Hellman problem (SIDH) on elliptic curves, was first published only in 2006, although the idea was kicking around a decade or so earlier. It combines a certain type of polynomial equation called an elliptic curve with a special type of graph. For the hard problem you are given two points in the graph and the goal is to find a path that goes from one to the other. It turns out that certain elliptic curves translate to particularly complicated graphs with lots of dead ends and cul-de-sacs. A big advantage of this type of cipher is that a lot of current cryptography already uses elliptic curves, although in a way that's not secure against quantum computers. The big question is whether these ciphers have been around long enough for potential weaknesses to have been found. There is one NIST alternate candidate, SIKE, which is based on SIDH.

SIDH could be very accessible to undergraduates. Many cryptography teachers already cover elliptic curve cryptography, although SIDH would definitely take more time to explain than the "traditional" version. There are several good expository papers on SIDH. "A friendly introduction to Supersingular Isogeny Diffie-Hellman" (Urbanik 2017) does a pretty good job of explaining the basic idea to someone who doesn't know elliptic curves at all. "Supersingular isogeny key exchange for beginners" (Costello 2020) requires a little more sophistication but goes into more depth. (It also has links to more advanced surveys and lecture notes.) Most cryptography textbooks now have a decent (or better) section on elliptic curves, including all the ones I've mentioned here, but they stop short of discussing isogenies.

7. **Conclusion**

NIST's announcement of the second-round candidates was postponed by the shutdown of the federal government and its announcement of the third-round candidates was delayed by COVID-19. What about the final selections? NIST has made it clear that there is likely to be more than one "winner" and also possibly a number of backups which look promising but aren't yet ready to be government standards. NIST's current plan is to make an announcement in Spring of 2022 with one group of algorithms to be standardized and another group that will advance to a fourth round of evaluation. If I had to bet, I'd go with one or two winners from LWE lattices and one from code-based cryptography, with SIDH as a backup. Draft standards for public comment are scheduled to be released in 2022 or 2023, and the final standards will hopefully be ready by 2024. Unless, of course, something new happens to hold up the process.

**References**

Alagic, G., J. Alperin-Sheriff, D. Apon, D. Cooper, Q. Dang, J. Kelsey, Y.-K. Liu, et al. 2020. *Status Report on the Second Round of the NIST Post-Quantum Cryptography Standardization Process*. NIST Internal or Interagency Report (NISTIR) 8309. Gaithersburg, MD: National Institute of Standards and Technology. https://csrc.nist.gov/publications/detail/nistir/8309/final.

Bernstein, D.J., J. Buchmann, and E. Dahmen, eds. 2009. *Post-Quantum Cryptography*. Berlin, Heidelberg: Springer. http://link.springer.com/10.1007/978-3-540-88702-7.

Costello, C. 2020. Supersingular Isogeny Key Exchange for Beginners. In *Selected Areas in Cryptography – SAC 2019*, ed. K.G. Paterson and D. Stebila, 21–50. Cham, Switzerland: Springer International Publishing.

Hoffstein, J., J. Pipher, and J.H. Silverman. 2014. *An Introduction to Mathematical Cryptography*. 2nd ed. Undergraduate Texts in Mathematics. New York, NY: Springer. http://dx.doi.org/10.1007/978-1-4939-1711-2.

Knospe, H. 2019. *A Course in Cryptography*. Pure and Applied Undergraduate Texts. Providence, Rhode Island: American Mathematical Society. http://www.librarything.com/work/23767116/book/173944779.

Peikert, C. 2016. *A Decade of Lattice Cryptography*. Cryptology ePrint Archive: Report 2015/939. IACR. https://eprint.iacr.org/2015/939.

Stallings, W. 2019. *Cryptography and Network Security: Principles and Practice*. 8th ed. Hoboken, NJ: Pearson Education, Inc.

Stinson, D.R., and M. Paterson. 2018. *Cryptography: Theory and Practice*. 4th Edition. Boca Raton, FL: Chapman and Hall/CRC.

Trappe, W., and L.C. Washington. 2020. *Introduction to Cryptography with Coding Theory*. 3rd ed. Hoboken, NJ: Pearson.

Urbanik, D. 2017. *A Friendly Introduction to Supersingular Isogeny Diffie-Hellman*. https://csclub.uwaterloo.ca/~dburbani/work/friendlysidh.pdf.